\documentclass [11pt] {article}
\usepackage{amssymb}
\usepackage{latexsym}
\usepackage{amsmath}
\title{Invariant Star Products of Wick Type:\\Classification and
Quantum Momentum Mappings}
\author{{\bf Michael F.\ M\"uller-Bahns
\thanks{bahns@mathi.uni-heidelberg.de}}\\[3mm]
Mathematisches Institut\\ Universit\"at
Heidelberg\\Im Neuenheimer Feld 288\\ D-69120 Heidelberg\\
Germany\\[3mm] {\bf Nikolai Neumaier
\thanks{neumaier@math.uni-frankfurt.de}}\\[3mm]
Fachbereich Mathematik\\ Universit\"at Frankfurt
\\Robert-Mayer-Stra\ss e 10\\D-60054 Frankfurt a.~M.\\
Germany\\[3mm]}
\date{Revised Version\\[3mm]February 2004}
\newcommand {\cc} [1] {\overline {{#1}}}
\newcommand {\Cinf} [1] {\mathcal C^\infty ({#1})}
\newcommand {\Ginf} [1] {\Gamma^\infty ({#1})}
\newcommand {\tr} [1]{\mathsf{tr} \left({#1}\right)}
\newcommand {\JN} [1] {J_0({#1})}
\newcommand {\JP} [1] {J_+({#1})}
\newcommand {\Jbold} [1] {J({#1})}
\newcommand {\im} {{\mathrm{i}}}

\newcommand {\Lie} {{\mathcal L}}
\newcommand {\ad} {{\mathrm{ad}}}
\newcommand {\Ad} {{\mathrm{Ad}}}
\newcommand{\starBT} {\mathbin{\star_{{\mbox{\rm \tiny BT}}}}}
\renewcommand {\d} {\mathrm{d}}

\newenvironment {PROOF}{\small{\sc Proof:}}{{\hspace*{\fill}
                       $\square$}}
\newtheorem {LEMMA} {Lemma} [section]

\newtheorem {PROPOSITION} [LEMMA] {Proposition}
\newtheorem {THEOREM} [LEMMA] {Theorem}
\newtheorem {COROLLARY} [LEMMA] {Corollary}

\newtheorem {REMARK}[LEMMA] {Remark}
\newtheorem {DEDUCTION}[LEMMA] {Deduction}
\newtheorem {EXAMPLE} [LEMMA] {Example}
\begin{document}
\maketitle
\begin{abstract}
\noindent We extend our investigations on $\mathfrak g$-invariant
Fedosov star products and quantum momentum mappings \cite{MN03a} to
star products of Wick type on pseudo-K\"ahler manifolds. Star
products of Wick type can be completely characterized by a local
description as given by Karabegov in \cite{Kar96} for star products
with separation of variables. We separately treat the action of a
Lie group $G$ on $\Cinf{M}[[\nu]]$ by (pull-backs with)
diffeomorphisms and the action of a Lie algebra $\mathfrak g$ on
$\Cinf{M}[[\nu]]$ by (Lie derivatives with respect to) vector
fields. Within Karabegov's framework, we prove necessary and
sufficient conditions for a given star product of Wick type to be
invariant in the respective sense. Moreover, our results yield a
complete classification of invariant star products of Wick type. We
also prove a necessary and sufficient condition for (the Lie
derivative with respect to) a vector field to be even a quasi-inner
derivation of a given star product of Wick type. We then transfer
our former results about quantum momentum mappings for $\mathfrak
g$-invariant Fedosov star products to the case of invariant star
products of Wick type.
\end{abstract}
\clearpage
\tableofcontents
\section{Introduction}
\label{IntroSec}
In all the existing approaches to quantization, the incorporation
of classical symmetries is a central issue that has proven to pose
serious problems. In the framework of deformation quantization,
however, this incorporation can at least be formulated very
naturally as it has already been indicated in the pioneering
articles \cite{BayFla78} by Bayen, Flato, Fr\o nsdal, Lichnerowicz,
and Sternheimer. Various notions of invariance of star products
with respect to actions of Lie groups and Lie algebras were
introduced and discussed by Arnal, Cortet, Molin, and Pinczon in
\cite{ArnCor83}. Previously, the existence of $G$-invariant
symplectic connections has been related to that of certain
$G$-invariant star products by Lichnerowicz in \cite{Lic80}. In
\cite{BerBieGut98} Bertelson, Bieliavsky, and Gutt proved that the
$G$-equivalence classes of $G$-invariant star products on a
symplectic manifold that possesses a $G$-invariant symplectic
connection are in bijection to formal series with values in the
second $G$-invariant de Rham cohomology of the manifold.

Another important notion in deformation quantization
-- which is one of the key ingredients of the formulation of phase
space reduction in this framework (cf.\ \cite{BorHerWal00,Fed98})
-- is that of a quantum momentum mapping, an object introduced and
studied in detail by Xu \cite{Xu98}. In \cite{MN03a}, we studied
$\mathfrak g$-invariant Fedosov star products and quantum momentum
mappings and we gave necessary and sufficient conditions for the
existence of quantum momentum mappings for these star products
which particularly showed that generally, the existence of a
classical momentum mapping does not imply the existence of a
quantum momentum mapping. (Some of the statements of \cite{MN03a}
have also been presented by Gutt in \cite{Gut02}, see also
\cite{GutRaw03}.)

The present letter extends these results to star products of Wick
type on pseudo-K\"ahler manifolds. Those are an important example
of star products compatible with an additional geometric structure,
namely the complex structure. They are closely related to geometric
quantization \cite{Woo91} with a complex polarization, and to
Berezin's quantization on K\"ahler manifolds. The star products
constructed in that context by Moreno \cite{Mor86} and by Cahen,
Gutt, and Rawnsley \cite{CahGutRaw93} are concrete examples.
Moreover, star products of Wick type also appear as asymptotic
expansions of the Berezin-Toeplitz quantization, cf.\
\cite{BorMeiSch94}, and for details \cite{KarSch00,Sch99}.

We begin by recalling the definition of these special star
products. On a pseudo-K\"ahler manifold $(M,\omega,I)$ a star
product is said to be of Wick type if the bidifferential operators
determining the star product contain only derivatives in
holomorphic directions in the first argument and only derivatives
in anti-holomorphic directions in the second argument. This
definition of star products of Wick type clearly is equivalent to
the condition that on every open subset $U\subseteq M$
star-right-multiplication with functions that are holomorphic on
$U$ and star-left-multiplication with functions that are
anti-holomorphic on $U$ both coincide with pointwise
multiplication. Star products of Wick type for general K\"ahler
manifolds were independently constructed by Karabegov in his work
on star products with separation of variables \cite{Kar96,Kar98,
Kar00} and by Bordemann and Waldmann \cite{BorWal97} using a
modified Fedosov construction with a fibrewise Wick product. The
latter results have been generalized in \cite{Neu02a}, where it is
shown under which conditions a generalized Fedosov star product
constructed using a fibrewise Wick product is of Wick type and that
in fact all star products of Wick type are generalized Fedosov star
products.

Our results constitute the starting point for the investigation of
the interplay of deformation quantization of Wick type and phase
space reduction of pseudo-K\"ahler manifolds. The first step in
this direction which is done in this letter is to derive necessary
and sufficient conditions that permit the formulation of
`reduction' of a star product of Wick type (by reduction of a star
product we mean the definition of a star product on the reduced
phase which is induced by a star product on the initial phase
space) that in the hitherto existing approaches is bound to the
$G$-invariance of the original star product and to the existence of
a quantum momentum mapping (cf.\ Remark \ref{ReductRem}).

The present letter is organized as follows: In Section
\ref{KarSec}, we briefly review some of Karabegov's main results as
needed in the sequel. We also collect some notational conventions
and shorthands.

In Section \ref{ClasSec}, we use Karabegov's results to
independently treat the cases of invariance of a given star product
of Wick type on a pseudo-K\"ahler manifold with respect to the
action of Lie groups and Lie algebras. We give a complete
classification of invariant star products of Wick type which shows
that for star products of Wick type, an even stronger result holds
than the one obtained in \cite{BerBieGut98} for general invariant
star products on a symplectic manifold. Eventually, we give a
necessary and sufficient condition for a vector field $X$ to be a
quasi-inner derivation (see below for a definition).

With the results from Section \ref{ClasSec}, our statements about
quantum momentum mappings from \cite{MN03a} are then adapted to the
case of invariant star products of Wick type in Section
\ref{QmmSec}. Finally, we provide some concrete examples of
invariant star products of Wick type, where the existence of a
quantum momentum mapping is guaranteed. In particular, the star
product corresponding to the Berezin-Toeplitz quantization has this
remarkable property provided there is a classical momentum mapping.

\noindent {\bf Conventions:} Whenever we speak of K\"ahler
manifolds $(M,\omega,I)$, we include pseudo-K\"ahler manifolds
since positivity of the Hermitean metric $\omega(\cdot,I\cdot)$ is
nowhere required. By $\Cinf{M}$, we denote the complex-valued
smooth functions and similarly $\Ginf{TM}$ stands for the
complex-valued smooth vector fields et cetera. Slightly abusing
notation we denote elements of $\Ginf{TM}$ and the corresponding
derivations on $\Cinf{M}$ by the same symbol. We use Einstein's
summation convention in local expressions.

\section{Karabegov's Description and Characterization of
Star Products of Wick Type}
\label{KarSec}
Let $(M,\omega,I)$ be a K\"ahler manifold of real dimension $2n$
with symplectic form $\omega$ and complex structure $I$. In a local
holomorphic chart, we write $Z_k:=\partial_{z^k}$ and
$\cc{Z}_l:=\partial_{\cc{z}^l}$ for local basis vector fields of
type $(1,0)$ and of type $(0,1)$ that locally span the $+\im$ and
$-\im$ eigenspaces $TM^{1,0}$ and $TM^{0,1}$ of the complex
structure $I$. For vector fields $X\in \Ginf{TM}=\Ginf{TM^{1,0}}
\oplus \Ginf{TM^{0,1}}$, we sometimes write $X=\chi +
\underline{\chi}$ with $\chi \in \Ginf{TM^{1,0}}$ and
$\underline{\chi}\in \Ginf{TM^{0,1}}$.

A star product of Wick type on $(M,\omega,I)$ then is by definition
a star product given for all $a,b \in \Cinf M$ by $a\star
b=\sum_{r=0}^{\infty}\nu^rC_r(a,b)$, where, using some coefficient
functions $C_r^{K;\cc{L}}$, in local holomorphic coordinates for
$r\geq 1$ each bidifferential operator has the form
$C_r(a,b)=\sum_{K,\cc{L}}C_r^{K;\cc{L}}
\frac{\partial^{|K|}a}{\partial z^K}
\frac{\partial^{|\cc{L}|}b}{\partial\cc{z}^L}$.

We shall now use Karabegov's results from \cite{Kar96} about star
products with separation of variables to give a unique description
of all star products of Wick type. Note that the star
products with separation of variables in Karabegovs original works
\cite{Kar96,Kar98,Kar00} differ from the star products of Wick type
considered here by a sign in the Poisson bracket and by an
interchange of the r\^{o}les of holomorphic and anti-holomorphic
coordinates. In the sequel, we shall therefore adapt Karabegov's
results to our setting.

Let now $\star$ be a star product of Wick type on $(M,\omega, I)$.
Then (cf.\ \cite[Lemma 2]{Kar96}) there are formal functions $u_k,
\cc{v}_l \in \Cinf{U}[[\nu]]$ defined on an open, contractible
domain $U \subseteq M$ of a local holomorphic chart $(z,U)$ of $M$
such that
\begin{equation}
\label{LinkRechWickMulKarEq}
a \star  u_k = a u_k + \nu Z_k(a)\qquad\textrm{ and } \qquad
\cc{v}_l \star  a = \cc{v}_l a + \nu \cc{Z}_l(a)
\end{equation}
for all $a\in \Cinf{U}[[\nu]]$. We shall always reserve the symbols
$u_k$ and $\cc{v}_l$ to denote functions as in
(\ref{LinkRechWickMulKarEq}). Moreover, Karabegov considers locally
defined formal series of one-forms $\alpha,\beta\in
\Ginf{T^*U}[[\nu]]$ given by $\alpha:= - u_k \d z^k$ which is of
type $(1,0)$ and $\beta:= \cc{v}_l \d
\cc{z}^l$ which is of type $(0,1)$. As the
$\star$-right-multiplication with $u_k$ obviously commutes with the
$\star$-left-multiplication with $\cc{v}_l$, one in addition
obtains from (\ref{LinkRechWickMulKarEq}) that
$\cc{\partial}\alpha= \partial \beta$. One can show that this
procedure yields a formal series of closed two-forms of type
$(1,1)$ on $M$. In the following, the so-defined formal two-form
that can be associated to any star product $\star$ of Wick type
shall be denoted by $K(\star)$. It is referred to as Karabegov's
characterizing form of the star product $\star$. It is easy to see
from the very definition that $K(\star)
\in \omega + \nu Z^2_{\mbox{\rm\tiny dR}}(M,\mathbb
C)^{1,1}[[\nu]]$, where $\nu Z^2_{\mbox{\rm\tiny dR}}(M,\mathbb
C)^{1,1}[[\nu]]=\{\Omega \in
\nu \Ginf{\mbox{$\bigwedge$}^2T^*M}[[\nu]]\,|\, \d \Omega=0,
\pi^{1,1}\Omega =\Omega\}$, and hence by the
$\partial$-$\cc{\partial}$-Poincar\'{e} Lemmas, there exist local
formal functions $\varphi\in\Cinf{U}[[\nu]]$ on every contractible
domain $U$ of holomorphic coordinates such that $K(\star)|_U =
\partial\cc{\partial} \varphi$; $\varphi$ is called a formal
local K\"ahler potential of $K(\star)$. With such a formal local
K\"ahler potential, the equations
\begin{equation}
\label{KahlPotStMultEq} a \star Z_k(\varphi) = a Z_k(\varphi) +
\nu Z_k(a)\qquad\textrm{ and } \qquad
\cc{Z}_l(\varphi) \star a = \cc{Z}_l(\varphi) a + \nu \cc{Z}_l(a)
\end{equation}
hold for all $a\in \Cinf{U}[[\nu]]$. Together, one finds:
\begin{THEOREM}{\rm\bf (\cite[Thm.\ 1]{Kar96})}
\label{KStarSatz} Let $\star$ be a star product of Wick type
on a K\"ahler manifold $(M,\omega,I)$. Then $K(\star) \in \omega +
\nu Z^2_{\mbox{\rm\tiny dR}}(M,\mathbb C)^{1,1}[[\nu]]$ associates
a formal series of closed two-forms of type $(1,1)$ on $M$  --
which is a deformation of the K\"ahler form $\omega$ -- to this
star product. In case $\varphi \in \Cinf{U}[[\nu]]$ is a formal
local K\"ahler potential of $K(\star)$, the equations
(\ref{KahlPotStMultEq}) hold for all $a\in \Cinf{U}[[\nu]]$.
\end{THEOREM}

Conversely, in \cite[Sect.\ 4]{Kar96} Karabegov has shown that to
each form $K$ as in the preceding theorem, one can assign a
star product $\star$ of Wick type such that the characterizing form
$K(\star)$ of this star product actually coincides with this given
$K$. To this end, Karabegov has given an explicit construction of
such a star product extensively using local considerations, and he
proved:

\begin{THEOREM}{\rm\bf (\cite[Thm.\ 2]{Kar96})}
Karabegov's characterizing form induces a bijection
\begin{equation}\label{WickMapOBij}
\{
\textrm{ star products of Wick type on $(M,\omega,I)$ }
\}\ni \star \quad\mapsto \quad K(\star)-\omega  \in
\nu Z^2_{\mbox{\rm\tiny dR}}(M,\mathbb C)^{1,1}[[\nu]]
\end{equation}
between star products of Wick type on $(M,\omega,I)$ and formal
closed two-forms on $M$ of type $(1,1)$ with formal degree $\geq
1$.
\end{THEOREM}

In \cite[Thm.\ 5.2]{Neu02a} an alternative proof of the fact that a
star product of Wick type is completely determined by one of the
equations (\ref{KahlPotStMultEq}) valid in every holomorphic chart
and hence that it is determined by its characterizing form has been
given.

\section{Classification of Invariant Star Products of Wick Type}
\label{ClasSec}

We are going to treat the cases of invariance with respect to group
actions on $\Cinf{M}[[\nu]]$ by diffeomorphisms and invariance with
respect to Lie algebra actions on $\Cinf{M}[[\nu]]$ by vector
fields separately and independently. That way, we do not have to
assume that the groups considered are connected (which would permit
to deduce the statement about the group action from the analogous
statement about the corresponding Lie algebra action). On the other
hand, we have to consider actions of Lie algebras anyway to derive
conditions for the existence of quantum momentum mappings, and it
turns out that our proofs also work for actions of complex Lie
algebras that, in general, trivially cannot be the infinitesimals
of actions of Lie groups. Nevertheless, the result about the Lie
algebra actions that are not infinitesimals of group actions should
rather be seen as a by-product since the interesting situation for
the application in phase space reduction is that of the action of a
Lie group.

Note that, in the case of a Lie group $G$, the statements for the
sufficient conditions for an ordinary Fedosov star product as
defined in \cite{Fed94} to be $G$-invariant are well-known (cf.\
\cite[Sect.\ 3]{BerBieGut98}). In \cite{MN03a}, we proved a
statement giving both necessary and sufficient conditions for
generalized Fedosov star products to be $\mathfrak g$-invariant.
Also note that, owing to the fact that every star product of Wick
type can be obtained using the generalized Fedosov construction on
K\"ahler manifolds \cite{Neu02a}, one could as well use a Fedosov
construction to prove the statements of the following sections. In
particular, the proofs of the sufficient conditions for invariance
with respect to a diffeomorphism and invariance with respect to a
vector field in Propositions \ref{DifAutKarProp} and
\ref{DerKarProp} are straightforward using Fedosov methods. In the
present letter, however, we avoid the use of the Fedosov
construction, and we give much less technical proofs only using
Karabegov's powerful description of star products of Wick type. The
methods we use are similar in spirit to those applied to flag
manifolds in \cite{Kar99}.

First we give necessary and sufficient conditions for a
diffeomorphism of $M$ to induce -- via pull-back -- an automorphism
of a given star product of Wick type.

\begin{PROPOSITION}
\label{DifAutKarProp}
Let $\star$ be a star product of Wick type on $(M,\omega,I)$ with
Karabegov's characterizing form $K(\star)$, and let $\phi$ be a
diffeomorphism of $M$. Then the pull-back $\phi^*$ is an
automorphism of $\star$ if and only if
\begin{equation}\label{DiffInvEq}
\phi^* I = I \quad \textrm{ and } \quad \phi^* K(\star)= K(\star).
\end{equation}
\end{PROPOSITION}
\begin{PROOF}
Let us first prove that $\phi^* I = I$ and $\phi^* K(\star)=
K(\star)$ are necessary conditions for $\phi^*$ to be an
automorphism of $\star$. Then ${\phi^{-1}}^*$ also being an
automorphism of $\star$ and Karabegov's relations
(\ref{LinkRechWickMulKarEq}) imply that on $\phi(U)$ we have
\begin{equation}
\label{phimininvEq}
a\star {\phi^{-1}}^*u_k={\phi^{-1}}^*(\phi^*a\star
u_k)={\phi^{-1}}^*\left(\phi^*a u_k+ \nu Z_k(\phi^* a)
\right)=a{\phi^{-1}}^*u_k + \nu \left({\phi^{-1}}^*Z_k\right)(a)
\end{equation}
for all $a \in \Cinf{\phi(U)}$. In particular, we can choose a
local holomorphic chart $(z',\phi(U))$, and then for $a=\cc{z'}^l$
the last equation becomes $\left({\phi^{-1}}^*Z_k
\right)(\cc{z'}^l)=0$, which implies that ${\phi^{-1}}^*Z_k$ is
still of type $(1,0)$, that is $I({\phi^{-1}}^*Z_k)=\im
{\phi^{-1}}^*Z_k$. But this implies that ${\phi^{-1}}^*$ maps
vector fields if type $(1,0)$ to vector fields of type $(1,0)$ and
hence $\phi^*I=I$. In particular, $\phi^*$ maps local holomorphic
charts to local holomorphic charts. Now we calculate $K(\star)$ and
$\phi^*K(\star)$. Using equation (\ref{phimininvEq}) for
$\phi^{-1}$ instead of $\phi$, we obtain $a \star \phi^*u_k
= a\phi^*u_k + \nu (\phi^*Z_k)(a)$ for $a \in \Cinf{\phi^{-1}(U)}$.
Since $\phi^* Z_k = Z'_k$ in the holomorphic chart $(z',U')=
\left(\phi^*z ,\phi^{-1}(U)\right)$, we can use $\phi^*u_k$ to
calculate $K(\star)$ in the chart $(z',U')$:
\begin{eqnarray*}
K(\star)|_{U'} & = & - \cc{\partial}(\phi^*u_k \d {z'}^k)
     =  - \cc{Z'_l}(\phi^*u_k) \d\cc{z'}^l\wedge\d {z'}^k
     =  - \phi^*(\cc{Z}_l(u_k)) \phi^* \d\cc{z}^l \wedge
        \phi^* \d z^k\\ & = & \phi^*(- \cc{\partial}(u_k \d z^k))
     =  \phi^*(K(\star)|_U).
\end{eqnarray*}
Analogously, one could use the equation $\phi^*\cc{v}_l \star b =
\phi^*\cc{v}_l b + \nu (\phi^*\cc{Z}_l)(b)$ for $b \in
\Cinf{\phi^{-1}(U)}$ to obtain the same result.

Conversely, let us now assume that $\phi^*I=I$ and that
$\phi^*K(\star)=K(\star)$. We want to prove that this implies that
$\phi^*(a \star b)=\phi^*a\star\phi^*b$ for all formal functions
$a,b \in \Cinf{M}[[\nu]]$. For this purpose, consider the star
product
\[
a\star' b :={\phi^{-1}}^*\left(\phi^*a\star \phi^*b\right).
\]
Since ${\phi^{-1}}^*$ is type-preserving and since ${\phi^{-1}}^*
K(\star) = K(\star)$, which in particular implies
${\phi^{-1}}^*\omega = \omega$, we obviously find that $\star'$ is
a star product of Wick type on $(M,\omega,I)$. But star products of
Wick type are uniquely determined by their Karabegov form $K$.
Therefore, the proof is done if we can show that
$K(\star')=K(\star)$, since then $\star'$ equals
$\star$ and hence $\phi^*$ is an automorphism of $\star$. To
see that this is true, we write ${\phi^{-1}}^* Z_k =
\widetilde{Z}_k$ in the holomorphic chart
$(\widetilde{z},\widetilde{U})=({\phi^{-1}}^*z,\phi(U))$ and
compute as in equation (\ref{phimininvEq}): $a \star'
{\phi^{-1}}^*u_k= a {\phi^{-1}}^*u_k + \nu \widetilde{Z}_k(a)$ for
$a\in
\Cinf{\phi(U)}$. Therefore, we have
\[
K(\star')|_{\phi(U)}= -
\cc{\widetilde{Z}}_l\left({\phi^{-1}}^*u_k\right)
\d\cc{\widetilde{z}}^l \wedge \d\widetilde{z}^k =
{\phi^{-1}}^*\left(-\cc{Z}_l(u_k) \d\cc{z}^l\wedge\d z^k \right) =
{\phi^{-1}}^*(K(\star)|_U)=K(\star)|_{\phi(U)},
\]
where we used $K(\star)= \phi^*K(\star)$ to obtain the last
equality.
\end{PROOF}

Now we derive necessary and sufficient conditions for a vector
field on $M$ to define -- via the Lie derivative -- a derivation of
a given star product of Wick type. The proof also gives an
important step of the proof for necessary and sufficient conditions
for a given derivation to be quasi-inner (cf.\ Proposition
\ref{quainnDerProp}).

\begin{PROPOSITION}
\label{DerKarProp}
Let $\star$ be a star product of Wick type on $(M,\omega,I)$ with
Karabegov's characterizing form $K(\star)$ and let $X  \in
\Ginf{TM}$ be a vector field on $M$. Then $X$ is a derivation of
$\star$ if and only if
\begin{equation}\label{LieInvEq}
\Lie_X I = 0
\quad\textrm{ and } \quad \Lie_X K(\star)=0.
\end{equation}
\end{PROPOSITION}
\begin{PROOF}
Let us first show that $\Lie_X I = 0 = \Lie_X K(\star)$ implies
that $X$ is a derivation of $\star$. We write $X=\chi +
\underline{\chi}$ with $\chi \in \Ginf{TM^{1,0}}$ and
$\underline{\chi}\in \Ginf{TM^{0,1}}$, locally $X=\chi^n Z_n +
\underline\chi^{\cc{l}} \cc{Z}_l$. Then we get from
$0=(\Lie_X I) Z_k$ that $\Lie_X Z_k$ must be of type $(1,0)$. Since
locally $\Lie_X Z_k=-\left(Z_k(\chi^n)Z_n +
Z_k(\underline\chi^{\cc{l}})\cc{Z}_l\right)$, this implies that
$Z_k(\underline\chi^{\cc{l}})=0$, which means that
$\underline\chi^{\cc{l}}$ is locally anti-holomorphic. Likewise,
from $(\Lie_X I)\cc{Z}_m=0$ we find that $\chi^n$ is locally
holomorphic. Now, since $K(\star)$ is closed, we have $0=\Lie_X
K(\star)=\d i_X K(\star)$, i.e.\ $i_X K(\star)$ is a closed
one-form on $M$. Hence on an open, contractible domain $U$ of a
chart, there exists a local formal function $a \in
\Cinf{U}[[\nu]]$ with $i_X K(\star)|_U=\d a$, hence $\partial a =
i_{\underline{\chi}}K(\star)|_U$ and
$\cc{\partial}a=i_{\chi}K(\star)|_U$. Using the first of these
equations with $K(\star)$ written out in local coordinates
$K(\star)|_U = Z_k(\cc{v}_l) \d z^k
\wedge \d \cc{z}^l$ then gives $Z_k(a)\d z^k = -
Z_k(\underline\chi^{\cc{l}}\cc{v}_l)\d z^k$, which follows since
$\underline\chi^{\cc{l}}$ is locally anti-holomorphic. Therefore,
locally $a=-\underline\chi^{\cc{l}}\cc{v}_l + \cc{h}$ with a
locally anti-holomorphic formal function $\cc{h}\in
\Cinf{U}[[\nu]]$. Then since $\star$ is of Wick type and since
$\underline\chi^{\cc{l}}$ is locally anti-holomorphic, we have
$\underline\chi^{\cc{l}}\cc{v}_l=\underline\chi^{\cc{l}}\star\cc{v}_l$.
Hence for all local functions $b \in \Cinf{U}$, the second of the
equations (\ref{LinkRechWickMulKarEq}) implies that the star
product $a\star b$ is given by
\[
a \star b=(\cc{h} - \underline\chi^{\cc{l}} \star \cc{v}_l) \star
b=\cc{h}b-
\underline\chi^{\cc{l}}(\cc{v}_l b + \nu\cc{Z}_l(b)).
\]
On the other hand, writing $K(\star)|_U=
-\cc{Z}_l (u_k) \d \cc{z}^l \wedge \d z^k$ and using
$\cc{\partial} a = i_{\chi}K(\star)|_U$, an analogous calculation
implies that $a=\chi^k u_k + h$ with a locally holomorphic formal
function $h\in \Cinf{U}[[\nu]]$. Then, as above, for all $b \in
\Cinf{U}$ we have
\[
b\star a=bh+ (bu_k+ \nu Z_k(b))\chi^k.
\]
Combining these equations, we find
\[
\ad_\star(a)b=\cc{h}b-
\underline\chi^{\cc{l}}(\cc{v}_l b + \nu\cc{Z}_l(b)) - bh - (bu_k+ \nu
Z_k(b))\chi^k.
\]
But evaluating this equation for $b=1$, we get $0=\cc{h} -
\underline\chi^{\cc{l}}{\cc{v}_l} - \chi^ku_k -h$, and therefore
$\ad_\star(a)b= - \nu\left(\underline\chi^{\cc{l}}\cc{Z}_l +
\chi^kZ_k
\right) (b) = -\nu X(b)$. But this means there is a local formal
function $a \in \Cinf{U}[[\nu]]$ such that on $U$ the Lie
derivative is quasi-inner, that is
\begin{equation}
\label{Lielocqi}
X|_{\Cinf{U}} = - \frac{1}{\nu}
\ad_\star(a)|_{\Cinf{U}}.
\end{equation}
We shall later use this fact to prove a necessary and sufficient
condition for $X$ to be globally quasi-inner. But now we have
finished the first part of the proof since (\ref{Lielocqi}) implies
that $X$ is a derivation of $\star$.

Conversely, let $X$ be a derivation of $\star$. Then, applying $X$
to the right-hand side of the first equation in
(\ref{LinkRechWickMulKarEq}), we have $X(a\star u_k)=X(a u_k) + \nu
X(Z_k(a))$, but using the derivation property and the first
equation of (\ref{LinkRechWickMulKarEq}) we also have $X(a\star
u_k)=X(a)u_k +
\nu Z_k (X(a)) + a\star X(u_k)$ and hence $a\star X(u_k)=a X(u_k) +
\nu [X,Z_k](a)$. Now in a local holomorphic chart $(z,U)$ of $M$,
we choose $a=\cc{z}^l$ and we write $X$ in local holomorphic
coordinates as $X=\chi^m Z_m+\underline\chi^{\cc{l}}\cc{Z}_l$.
First, $\star$-left-multiplication with $a=\cc{z}^l$ is just the
pointwise product, hence we get $0=\nu[X,
Z_k](\cc{z}^l)=-Z_k(\underline\chi^{\cc{l}})$. Analogously, from
the second equation of (\ref{LinkRechWickMulKarEq}) one finds
$-\cc{Z}_l(\chi^m)=0$. Now we have in local holomorphic
coordinates: $(\Lie_X I)Z_k
= \im \Lie_X Z_k  - I([\chi^m Z_m+\underline\chi^{\cc{l}}
\cc{Z}_l,Z_k])$. But $Z_k(\underline\chi^{\cc{l}})=0$ implies
$[\chi^m Z_m+\underline\chi^{\cc{l}}\cc{Z}_l,Z_k]=-Z_k(\chi^m)Z_m$.
Therefore $I\Lie_X Z_k= \im \Lie_X Z_k$ and hence $(\Lie_XI)Z_k=0$.
Likewise, one finds $(\Lie_XI)\cc{Z}_l=0$. Together this implies
$\Lie_X I=0$. To show that $\Lie_X K(\star)=0$, we use that
$K(\star)=
\cc{\partial}(-u_k\d z^k)$ and that we have $\Lie_X I =0$. Then
$\Lie_X K(\star)= -\cc{\partial}\Lie_X(u_k\d
z^k)=-\cc{\partial}(X(u_k)\d z^k + u_k \d \chi^k)=-(\cc{Z}_l
(X(u_k))
\d\cc{z}^l\wedge \d z^k + \cc{Z}_l(u_k) \d\cc{z}^l\wedge \d
\chi^k + u_k \cc{\partial}\d \chi^k)= -(\cc{Z}_l(X(u_k)) \d
\cc{z}^l\wedge \d z^k + \cc{Z}_l(u_m)Z_k(\chi^m)\d \cc{z}^l\wedge
\d z^k)$, where we have used that $\cc{\partial}\d \chi^k=0$ since
$\chi^k$ is locally holomorphic. But this last expression equals
zero, since applying $X$ to $\cc{v}_l\star u_k$, again the
derivation property and (\ref{LinkRechWickMulKarEq}) yield $[X,
Z_k](\cc{v}_l)=\cc{Z}_l(X(u_k))$. Explicitly this becomes
$-Z_k(\chi^m)\cc{Z}_l(u_m)=\cc{Z}_l(X(u_k))$ since $[X,
Z_k]=-Z_k(\chi^m)Z_m$, and $Z_m(\cc{v}_l)=\cc{Z}_l(u_m)$. Hence
$\Lie_X K(\star)=0$, and the proposition is proved.
\end{PROOF}

We have now proved all the prerequisites we need to give a complete
classification of all the invariant star products of Wick type on a
K\"ahler manifold. We first have to recall some definitions of
notions of invariance of star products from \cite{ArnCor83}.

Let $G$ be a Lie group and let $\Phi : G \times M \to M$ denote a
(left-)action of $G$ on $M$. Writing $\phi_g$ $\forall g \in G$ for
the diffeomorphism of $M$ defined by $\phi_g(m):=\Phi(g,m)$
$\forall m \in M$, obviously $r(g)a :=
\phi_{g^{-1}}^*a$ defines a Lie group action of $G$ on $\Cinf{M}$
that naturally extends to a Lie group action on $\Cinf{M}[[\nu]]$.
Now recall that a star product $\star$ on $(M,\omega)$ is called
$G$-invariant if $r(g)$ is an automorphism of $\star$ for all $g\in
G$.

Furthermore, let $\mathfrak g$ be a finite dimensional real or complex
Lie algebra and let $X_\cdot : \mathfrak g \to \Ginf{TM}: \xi
\mapsto X_\xi$ denote a Lie algebra anti-homomorphism. Then
obviously $\varrho(\xi)a:=- X_\xi (a)$ defines a Lie algebra action
of $\mathfrak g$ on $\Cinf{M}$ that also extends naturally to
$\Cinf{M}[[\nu]]$. Also recall that a star product $\star$ on
$(M,\omega)$ is called $\mathfrak g$-invariant if $\varrho(\xi)$ is
a derivation of $\star$ for all $\xi \in \mathfrak g$.

Observe that from the action of a Lie group $G$ one obtains a
corresponding Lie algebra action of $\mathfrak g = \mathrm{Lie}
(G)$ by $X_\xi(m): = \left.\frac{d}{dt}\right|_{t=0}\Phi
(\exp(t\xi),m)$ for all $m\in M$ and in this case $G$-invariance of
a star product clearly implies $\mathfrak g$-invariance. In case
$G$ is in addition connected one even has that $\mathfrak
g$-invariance with respect to the above action implies
$G$-invariance.

With the notations from above and from Propositions
\ref{DifAutKarProp} and \ref{DerKarProp} we obtain the following
classification result:

\begin{THEOREM}
\label{ClasInvKarThm}
Let $(M,\omega,I)$ be a K\"ahler manifold.
\begin{enumerate}
\item
For a given Lie group action $r$ as above there are $G$-invariant
star products of Wick type on $(M,\omega,I)$ if and only if
$\phi_g^*I = I$ and $\phi_g^*\omega = \omega$ for all $g\in G$. In
this case the set of $G$-invariant star products of Wick type is in
bijection to $\{
\Omega  \in \nu Z^2_{\mbox{\rm\tiny dR}}(M,\mathbb C)^{1,1}[[\nu]]
\,|\, \phi_g^* \Omega = \Omega\,\, \forall g \in G\}$ and the
bijection is given by the restriction of the mapping according to
(\ref{WickMapOBij}) to the respective $G$-invariant elements.
\item
For a given Lie algebra action $\varrho$ as above there are
$\mathfrak g$-invariant star products of Wick type on
$(M,\omega,I)$ if and only if $\Lie_{X_\xi}I = 0$ and
$\Lie_{X_\xi}\omega = 0$ for all $\xi \in \mathfrak g$. In this
case the set of $\mathfrak g$-invariant star products of Wick type
is in bijection to $\{
\Omega  \in \nu Z^2_{\mbox{\rm\tiny dR}}(M,\mathbb C)^{1,1}[[\nu]]
\,|\,  \Lie_{X_\xi}\Omega = 0\,\,\forall \xi \in \mathfrak g\}$
and the bijection is given by the restriction of the mapping
according to (\ref{WickMapOBij}) to the respective $\mathfrak
g$-invariant elements.
\end{enumerate}
\end{THEOREM}

\begin{REMARK}\label{ReductRem}
It is interesting to note that the necessary and sufficient
conditions for the existence of $G$-invariant star products of Wick
type on K\"ahler manifolds given above are precisely -- besides
other additional conditions that guarantee that the reduced phase
space exists as a smooth symplectic manifold or more generally as a
stratified symplectic space -- the conditions that are sufficient
for the reduced phase space to be a K\"ahler manifold or more
generally a stratified K\"ahler space (cf.\ \cite{HeiHuc00}). In
particular, this means that the notion of star products of Wick
type can also be defined on the reduced phase space and therefore
the question whether `deformation quantization of Wick type
commutes with reduction' can at least be given a meaning. Note,
however, that unless there is a quantum momentum mapping, one does
not even have the possibility to say what the relation between a
star product on the reduced phase space and a reduced star product
should be, since so far there is no obvious method to obtain the
latter in that case.
\end{REMARK}

\noindent We now state a proposition giving a necessary and
sufficient condition for a derivation of a star product of Wick
type given by a vector field $X$ to satisfy
\begin{equation}\label{quainnDefEq}
X(b)=-\frac{1}{\nu}\ad_\star(a)b
\end{equation}
with some $a \in \Cinf{M}[[\nu]]$ for all $b\in \Cinf{M}[[\nu]]$.
In this case $X$ is called a quasi-inner (or essentially inner
\cite{GutRaw99}) derivation.

The condition given in \cite{MN03a} was first presented for
ordinary slightly more special Fedosov star products by Gutt in
\cite{Gut02} (cf.\ also \cite[Thm.\ 7.2]{GutRaw03}). The proof that
the given condition is sufficient for these Fedosov star products
was originally published in \cite[Prop.\ 4.3]{Kra98}.

\begin{PROPOSITION}{\rm\bf (cf.\
{\cite[Prop.\ 3.9]{MN03a}})} \label{quainnDerProp} Let $\star$ be a
star product of Wick type on $(M, \omega, I)$ with Karabegov's
characterizing form $K(\star)$. Assume that $X\in \Ginf{TM}$ is a
vector field such that $X$ is a derivation of $\star$. Then $X$ is
even quasi-inner if and only if there is a formal function $a
\in \Cinf{M}[[\nu]]$ such that
\begin{equation}\label{qinnEq}
\d a = i_X K(\star),
\end{equation}
and then $X(b) =X^\omega_{a_0}(b)= - \frac{1}{\nu}
\ad_\star(a)b$ for all $b\in \Cinf{M}[[\nu]]$, where we have
written $a = a_0 + a_+$ with $a_0\in \Cinf{M}$ and $a_+ \in \nu
\Cinf{M}[[\nu]]$, and $X^\omega_{a_0}$ denotes the Hamiltonian
vector field of $a_0$ with respect to $\omega$.
\end{PROPOSITION}
\begin{PROOF}
Let $X$ be a derivation of $\star$, then from Proposition
\ref{DerKarProp} we know that $\Lie_X I = 0
= \Lie_X K(\star)$. Now we have already shown in equation
(\ref{Lielocqi}) that under these conditions on every open,
contractible set $U \subseteq M$ there is a local formal function
$a \in \Cinf{U}[[\nu]]$ with $\d a=i_X K(\star)|_U$ and that this
implies that $X|_{\Cinf{U}}= - \frac{1}{\nu}
\ad_\star(a)|_{\Cinf{U}}$. Hence if there is a globally defined
formal function $a \in \Cinf{M}[[\nu]]$ with $\d a=i_XK(\star)$,
the proof of Proposition \ref{DerKarProp} also shows that $X$ is
quasi-inner.

Conversely let $X=\chi+\underline{\chi}$, where $\chi
\in \Ginf{TM^{1,0}}$ and $\underline{\chi} \in
\Ginf{TM^{0,1}}$, be a quasi-inner derivation of $\star$, that
is for all $b \in \Cinf{M}$ we have $X (b) =
-\frac{1}{\nu} \ad_\star(a)b$, or equivalently
$b\star a - \nu \chi(b)= a\star b + \nu \underline{\chi}(b)$. Since
$\star$ is of Wick type, the left-hand side of the former equation
contains only derivatives of $b$ in holomorphic directions, while
the right-hand side only contains derivatives in anti-holomorphic
ones. This implies that $b\star a - \nu \chi(b)= c b  (=a\star b +
\nu \underline{\chi}(b))$, with a formal function $c
\in \Cinf{M}[[\nu]]$ and for $b=1$ this particularly yields $c=a$
and therefore we have $\chi(b)=\frac{1}{\nu}( b \star a- b a)$ and
$\underline{\chi}(b)=\frac{1}{\nu}(a b - a\star b)$. As in the
previous proofs, we now use the equations
(\ref{LinkRechWickMulKarEq}) and calculate $\underline{\chi}( u_k)
=\frac{1}{\nu}(a u_k - (a u_k + \nu Z_k(a)))= - Z_k(a)$ and likewise
$\chi(\cc{v}_l)=\cc{Z}_l(a)$. From these equations and the very
definition of $K(\star)$ it is obvious that $\partial a =
i_{\underline{\chi}}K(\star)$ and likewise $\cc{\partial}
a=i_{\chi}K(\star)$ and hence (\ref{qinnEq}) holds and the
necessary condition is proved. For the remaining statement of the
proposition, just observe that the zeroth order in $\nu$ implies
that $\d a_0 = i_X \omega$, and hence $X=X^\omega_{a_0}$ is the
Hamiltonian vector field of the function $a_0 \in \Cinf{M}$.
\end{PROOF}

\section{Quantum Momentum Mappings for Invariant Star Products
of Wick Type}\label{QmmSec}

We are now in the position to derive the analogues of our results
on quantum momentum mappings for Fedosov star products \cite[Sect.\
4]{MN03a}.

For the sake of brevity, and since all the proofs are (almost)
literally the same, we simply write down the statements themselves
and refer the reader to \cite[Sect.\ 4]{MN03a} for proofs and some
additional comments.

We no longer treat the cases of Lie group actions and Lie algebra
actions separately and we simply speak of invariant star products
in this section. Moreover, in the case of a Lie group action, we
always consider the induced Lie algebra action by means of its
infinitesimal generators.

First we shall need to recall some definitions and notations from
\cite{Xu98} and \cite{MN03a}. Considering some complex vector space
$V$ endowed with a representation $\pi: \mathfrak g \to
\mathsf{Hom}(V,V)$ of the Lie algebra $\mathfrak g$ in $V$, we
denote the space of $V$-valued $k$-multilinear alternating forms on
$\mathfrak g$ by $C^k(\mathfrak g,V)$, and the corresponding
Chevalley-Eilenberg differential is denoted by $\delta_\pi
: C^\bullet(\mathfrak g,V)\to C^{\bullet +1}(\mathfrak g,V)$.
Moreover, the spaces of the corresponding cocycles, coboundaries,
and the corresponding cohomology spaces shall be denoted by
$Z_\pi^k(\mathfrak g,V)$, $B_\pi^k(\mathfrak g,V)$, and
$H_\pi^k(\mathfrak g,V)$, respectively.

A Lie group action $r$ or a Lie algebra action $\varrho$ is called
Hamiltonian if and only if there is an element $J_0 \in
C^1(\mathfrak g,\Cinf{M})$ such that $X^\omega_{\JN{\xi}} = X_\xi$
for all $\xi \in \mathfrak g$, i.e.\ $i_{X_\xi} \omega= \d
\JN{\xi}$. In this case $\varrho(\xi)\cdot =
\{\JN{\xi},{}\cdot{}\}$, and $J_0$ is said to be a Hamiltonian for
the action $r$ resp. $\varrho$. In case $J_0$ is $\mathfrak
g$-equivariant, i.e.\ $\varrho(\xi)
\JN{\eta}=\{\JN{\xi},\JN{\eta}\}= \JN{[\xi,\eta]}$ for all
$\xi,\eta\in \mathfrak g$, the Hamiltonian $J_0$ is called a
classical momentum mapping. For a Lie group action, this is the case
if in particular we have $r(g)\JN{\xi}= \JN{\Ad(g)\xi}$ for all
$g\in G$ and all $\xi \in \mathfrak g$, and in this case $J_0$ is
called a $G$-equivariant classical momentum mapping.

Let now $\star$ be an invariant star product. Then $J = J_0 +
J_+\in C^1(\mathfrak g,\Cinf{M})[[\nu]]$ with $J_0\in C^1(\mathfrak
g,\Cinf{M})$ and $J_+\in\nu C^1(\mathfrak g,\Cinf{M})[[\nu]]$ is
called a quantum Hamiltonian for the action $r$ resp. $\varrho$ in
case
\begin{eqnarray}\label{preqmmEq}
\varrho (\xi) = \frac{1}{\nu} \ad_\star (\Jbold{\xi})\quad
\textrm{for all}\quad\xi \in \mathfrak g.
\end{eqnarray}
$J$ is called a quantum momentum mapping if in addition
\begin{eqnarray}\label{qmmEq}
\frac{1}{\nu}\left(\Jbold{\xi}\star \Jbold{\eta}-
\Jbold{\eta}\star \Jbold{\xi}\right) = \Jbold{[\xi,\eta]}
\end{eqnarray}
for all $\xi,\eta \in \mathfrak g$, i.e.\ in case $J$ is $\mathfrak
g$-equivariant.

\noindent The zeroth orders of (\ref{preqmmEq}) and (\ref{qmmEq})
just mean that $J_0$ is a classical momentum mapping. As for the
classical momentum mapping $J_0$ a quantum Hamiltonian $J$ for a
Lie group action $r$ that is $G$-equivariant, i.e.\ that satisfies
$r(g)J(\xi) = J(\Ad(g)\xi)$ for all $g\in G$ and for all $\xi \in
\mathfrak g$, clearly defines a quantum momentum mapping which is
called a $G$-equivariant quantum momentum mapping.

Also recall the definition of a strongly invariant star product
from \cite{ArnCor83}: Let $J_0$ be a classical momentum mapping for
the action $r$ resp. $\varrho$. Then an invariant star product is
called strongly invariant if and only if $J=J_0$ defines a quantum
Hamiltonian for this action.

\noindent As in \cite[Ded.\ 4.4]{MN03a}, the condition for the
existence of a quantum Hamiltonian for an action $r$ resp.
$\varrho$ can be seen directly from Proposition
\ref{quainnDerProp}:
\begin{DEDUCTION}\label{qHamDed}
An invariant star product $\star$ of Wick type on $(M,\omega,I)$
with Karabegov's characterizing form $K(\star)$ admits a quantum
Hamiltonian if and only if there is an element $J
\in C^1(\mathfrak g,\Cinf{M})[[\nu]]$ such that
\begin{equation}\label{preqmmKarEq}
\d \Jbold{\xi} = i_{X_\xi}K(\star)\quad \forall \xi \in
\mathfrak g,
\end{equation}
i.e.\ if and only if the cohomology class of $i_{X_\xi}K(\star)$
vanishes for all $\xi
\in \mathfrak g$. Moreover, from equation (\ref{preqmmKarEq}) $J$
is determined (in case it exists) up to elements in $C^1(\mathfrak
g,\mathbb C)[[\nu]]$.
\end{DEDUCTION}

\noindent We also have the following statement on strong invariance:
\begin{COROLLARY}{\rm\bf(cf.\ \cite[Cor.\ 4.7]{MN03a})}
Let $J_0$ be a classical momentum mapping for the action $r$ resp.
$\varrho$. Then an invariant star product of Wick type $\star$ with
Karabegov's characterizing form $K(\star)$ is strongly invariant if
and only if
\begin{equation}\label{iXOmNul}
i_{X_\xi}(K(\star)-\omega) = 0 \quad\textrm{for all}\quad\xi \in
\mathfrak g.
\end{equation}
In this case every classical momentum mapping defines a quantum
momentum mapping for $\star$.
\end{COROLLARY}

\noindent In the general case of an invariant star
product of Wick type, \cite[Prop.\ 4.8]{MN03a} and \cite[Cor.\
4.10]{MN03a} also transfer directly:

\begin{PROPOSITION}{\rm\bf (cf.\ \cite[Prop.\
4.8]{MN03a})}\label{QmmProp} Let $J$ be a quantum Hamiltonian for a
star product of Wick type $\star$ with Karabegov's characterizing
form $K(\star)$, then $\lambda \in C^2(\mathfrak
g,\Cinf{M})[[\nu]]$ defined by
\begin{equation}\label{lamDefEq}
\lambda(\xi,\eta):= \frac{1}{\nu}\left(\Jbold{\xi}\star
\Jbold{\eta}- \Jbold{\eta} \star \Jbold{\xi}\right) -
\Jbold{[\xi,\eta]}
\end{equation}
lies in $C^2(\mathfrak g,\mathbb C)[[\nu]]$ and is an element of
$Z_0^2(\mathfrak g,\mathbb C)[[\nu]]$ which is explicitly given by
\begin{equation}\label{lamexplEq}
\lambda(\xi,\eta)=K(\star)(X_\xi, X_\eta)
- \Jbold{[\xi,\eta]}.
\end{equation}
The cohomology class $[\lambda]\in H_0^2(\mathfrak g,\mathbb C
)[[\nu]]$ does not depend on the choice of $J$. Moreover, quantum
momentum mappings exist if and only if $[\lambda]=[0]\in
H_0^2(\mathfrak g,\mathbb C)[[\nu]]$, and for every $\tau\in
C^1(\mathfrak g,\mathbb C)[[\nu]]$ such that $\delta_0 \tau =
\lambda$ the element $J^\tau :=
J - \tau\in C^1(\mathfrak g,\Cinf{M})[[\nu]]$ is a quantum momentum
mapping for $\star$. Finally, the quantum momentum mapping (if it
exists) is unique up to elements in $Z_0^1(\mathfrak g,\mathbb
C)[[\nu]]$, and hence we have uniqueness if and only if
$H_0^1(\mathfrak g,\mathbb C)=0$.
\end{PROPOSITION}

\begin{REMARK}In the case of a $G$-action a quantum Hamiltonian $J$
clearly defines a smooth mapping $\sigma : G \to C^1(\mathfrak g,
\mathbb C)[[\nu]]$ by $(\sigma(g))(\xi) : = r(g)\Jbold{\xi} -
\Jbold{\Ad(g)\xi}$ which is a cocycle on $G$, i.e.\ $\sigma(g g')
= \sigma(g') + \sigma(g) \circ \Ad(g')$ for all $g,g'\in G$.
Moreover, the cohomology class $[\sigma]$ is independent of the
choice of $J$, i.e.\ for a different quantum Hamiltonian $J'$ we
have $\sigma'(g) - \sigma (g) = \gamma - \gamma \circ \Ad(g)$,
where $\gamma \in C^1(\mathfrak g,\mathbb C)[[\nu]]$. Clearly, a
$G$-equivariant quantum momentum mapping exists if and only if the
cohomology class of $\sigma$ vanishes, i.e.\ if and only if there
is an element $\tau \in C^1(\mathfrak g,\mathbb C)[[\nu]]$ such
that $\sigma(g) = \tau - \tau \circ \Ad (g)$ and in this case
$J^\tau := J - \tau$ is a $G$-equivariant quantum momentum mapping
which is unique up to $G$-invariant elements in $C^1(\mathfrak
g,\mathbb C)[[\nu]]$.
\end{REMARK}

\noindent We would like to recall the example of semi-simple Lie
algebras we gave in \cite[Ex.\ 4.9]{MN03a}. In this case, owing to
the Whitehead Lemmas and a deformed version of Sternberg's Lemma,
there is a unique quantum momentum mapping for every invariant star
product of Wick type. Another important example clearly is that of
Abelian Lie algebras, where, even if a quantum momentum mapping
exists, it is never unique. Yet another interesting example, where
the existence of a quantum momentum mapping is guaranteed and which
is specific for the case of K\"ahler manifolds considered in the
present letter is the following:

\begin{EXAMPLE}
In case there is a global formal
K\"ahler potential $\varphi \in
\Cinf{M}[[\nu]]$ for $K(\star)$, i.e.\ $K(\star)= \partial
\cc{\partial}\varphi$, that additionally satisfies $X_\xi(\varphi)=0$
for all $\xi\in\mathfrak g$, it is a straightforward verification
to see that $J$ defined by $\Jbold{\xi}:= \frac{1}{2\im}(I
X_\xi)(\varphi)$ fulfills $\d \Jbold{\xi}= i_{X_\xi}K(\star)$ and
$\Jbold{[\xi,\eta]} = K(\star)(X_\xi,X_\eta)$ and therefore defines
a quantum momentum mapping for the invariant star product $\star$.
If in the case of a group action $\varphi$ additionally satisfies
$r(g)\varphi = \varphi$ for all $g\in G$, then the quantum momentum
mapping $J$ defined above evidently is $G$-equivariant.
\end{EXAMPLE}

Finally, in the case of star products of Wick type we also find that in general
the existence of a classical momentum mapping does not imply the
existence of a quantum momentum mapping:

\begin{COROLLARY}{\rm\bf (cf.\ \cite[Cor.\ 4.10]{MN03a})}
Let $\star$ be an invariant star product of Wick type with
Karabe\-gov's characterizing form $K(\star)$, and assume that there
is a classical momentum mapping $J_0$ for the action $r$ resp.\ $\varrho$,
then a quantum momentum mapping $J$ exists if and only if there is
an element $J_+\in \nu C^1(\mathfrak g,\Cinf{M})[[\nu]]$ such that
\begin{equation}\label{qmmExClassExEq}
i_{X_\xi} (K(\star)-\omega) = \d \JP{\xi}\quad\textrm{ and }\quad
(K(\star)-\omega)(X_\xi,X_\eta) = (\delta_\varrho J_+)(\xi,\eta)
\quad \forall \xi,\eta \in \mathfrak g,
\end{equation}
and these equations determine $J_+$ up to elements of $\nu
Z_0^1(\mathfrak g,\mathbb C)[[\nu]]$.
\end{COROLLARY}

We want to conclude our considerations giving an exceptional
example of an invariant but in general not strongly invariant star
product of Wick type, where the existence of a classical momentum
mapping already guarantees the existence of a quantum momentum
mapping.

\begin{EXAMPLE}
In \cite{KarSch00} it has been shown that the star product
$\starBT$ of Wick type that arises from the asymptotic expansion of
the Berezin-Toeplitz quantization corresponds to the characterizing
form $K(\starBT) = \omega + \frac{2\nu}{\im} \rho$, where $\rho$
denotes the Ricci form which -- using the curvature $R$ of the
K\"ahler connection $\nabla$ -- is explicitly given by $\rho(Y,Y')
= - \frac{1}{4}\tr{R(Y,Y') I}$ for $Y,Y'\in \Ginf{TM}$. In case
$\omega$ and $I$ are invariant evidently the K\"ahler connection
and its curvature are invariant implying that $\starBT$ is an
invariant star product by Theorem \ref{ClasInvKarThm}. Using the
invariance of $\nabla$ in the form $[X_\xi, \nabla_Y Y'] =
\nabla_{[X_\xi,Y]} Y' + \nabla_Y [X_\xi,Y']$ together with the
property $\nabla I = 0$ of the K\"ahler connection and the
invariance of $I$ it is an easy computation to obtain that
$i_{X_\xi} \rho = \d j (\xi)$, where $j(\xi) = \frac{1}{4}
\mathsf{div} (IX_\xi)$ and $\mathsf{div}$ denotes the covariant
divergence with respect to $\nabla$, i.e.\ $\mathsf{div}(Y) =
\tr{\nabla Y}$ for $Y\in \Ginf{TM}$. Moreover, using the invariance
of $\nabla$ and $I$ once again it is rather obvious that $\rho
(X_\xi, X_\eta) = j ([\xi,\eta])$ and in the case of a group action
we even have $r(g)j (\xi) = j (\Ad(g)\xi)$. Consequently, in this
situation the existence of a ($G$-equivariant) classical momentum
mapping $J_0$ implies that $J
= J_0 + \frac{2\nu}{\im} j$ defines a ($G$-equivariant) quantum momentum
mapping.
\end{EXAMPLE}
\begin{small}

\end{small}

\begin{thebibliography}{99}
%
\bibitem{ArnCor83}
{\sc Arnal, D., Cortet, J.\ C., Molin, P., Pinczon, G.:} {\it
Covariance and geometrical invariance in $*$ quantization.} J.\
Math.\ Phys.\ {\bf 24}, 276--283 (1983).
%
\bibitem{BayFla78}
{\sc Bayen, F., Flato, M., Fr\o nsdal, C., Lichnerowicz, A.,
Sternheimer, D.:} {\it Deformation Theory and Quantization.} Ann.\
Phys.\ {\bf 111}, Part I: 61--110, Part II: 111--151 (1978).
%
\bibitem{BerBieGut98}
{\sc Bertelson, M., Bieliavsky, P., Gutt, S.:} {\it Parametrizing
Equivalence Classes of Invariant Star Products.} Lett.\ Math.\
Phys.\ {\bf 46}, 339--345 (1998).
%
\bibitem{BorMeiSch94}
{\sc Bordemann, M., Meinrenken, E., Schlichenmaier, M.:} {\it
Toeplitz quantization of K\"ahler manifolds and $\mathrm{gl}(N)$,
$N \to\infty$ limits.} Commun.\ Math.\ Phys.\ {\bf 165}, 281--296
(1994).
%
\bibitem{BorHerWal00}
{\sc Bordemann, M., Herbig, H.-C., Waldmann, S.:} {\it BRST
Cohomology and Phase Space Reduction in Deformation Quantization.}
Commun.\ Math.\ Phys.\ {\bf 210}, 107--144 (2000).
%
\bibitem{BorWal97}
{\sc Bordemann, M., Waldmann, S.:} {\it A Fedosov Star Product of
Wick Type for K\"ahler Manifolds.} Lett.\ Math.\ Phys.\ {\bf 41},
243--253 (1997).
%
\bibitem{CahGutRaw93}
{\sc Cahen, M., Gutt, S., Rawnsley, J.:} {\it Quantization of
K\"ahler manifolds.\ II.} Trans.\ Amer.\ Math.\ Soc.\ {\bf 337, 1},
73--98 (1993).
%
\bibitem{Fed94}
{\sc Fedosov, B.\ V.:} {\it A Simple Geometrical Construction of
Deformation Quantization.} J.\ Diff.\ Geom.\ {\bf 40}, 213--238
(1994).
%
\bibitem{Fed98}
{\sc Fedosov, B.\ V.:} {\it Non-Abelian Reduction in Deformation
Quantization.} Lett.\ Math.\ Phys.\ {\bf 43}, 137--154 (1998).
%
\bibitem{GutRaw99}
{\sc Gutt, S., Rawnsley, J.:} {\it Equivalence of star products on
a symplectic manifold; an introduction to Deligne's \v{C}ech
cohomology classes.} J.\ Geom.\ Phys.\ {\bf 29}, 347--392 (1999).
%
\bibitem{Gut02}
{\sc Gutt, S.:} {\it Star products and group actions.} Contribution
to the Bayrischzell Workshop, April 26--29, 2002.
%
\bibitem{GutRaw03}
{\sc Gutt, S., Rawnsley, J.:} {\it Natural star products on
symplectic manifolds and quantum moment maps.} Preprint, April
2003, {\bf math.SG/0304498 v1}, to appear in Lett.\ Math.\ Phys..
%
\bibitem{HeiHuc00}
{\sc Heinzner, P, Huckleberry, A.:} {\it K\"ahlerian structures on
symplectic reductions.} in: {\sc Peternell, T.\ (ed.):} Complex
analysis and algebraic geometry.\ A volume in memory of Michael
Schneider.\ Walter de Gruyter, Berlin, 225--253 (2000).
%
\bibitem{Kar96}
{\sc Karabegov, A.\ V.:} {\it Deformation Quantization with
Separation of Variables on a K\"ahler Manifold.} Commun.\ Math.\
Phys.\ {\bf 180}, 745--755 (1996).
%
\bibitem{Kar98}
{\sc Karabegov, A.\ V.:} {\it Cohomological Classification of
Deformation Quantization with Separation of Variables.} Lett.\
Math.\ Phys.\ {\bf 43}, 347--357 (1998).
%
\bibitem{Kar99}
{\sc Karabegov, A.\ V.:} {\it Pseudo-K\"ahler Quantization on Flag
Manifolds.} Commun.\ Math.\ Phys.\ {\bf 200}, 355--379 (1999).
%
\bibitem{Kar00}
{\sc Karabegov, A.\ V.:} {\it On Fedosov's approach to Deformation
Quantization with Separation of Variables.} in: {\sc Dito, G.,
Sternheimer, D.\ (eds.):} Conf\'{e}rence Mosh\'{e} Flato 1999,
Vol.\ II.\ Kluwer Academic Publ., Dordrecht, 167--176 (2000).
%
\bibitem{KarSch00}
{\sc Karabegov, A.\ V., Schlichenmaier, M.:} {\it Identification of
Berezin-Toeplitz Deformation Quantization.} J.\ reine angew.\
Math.\ {\bf 540}, 49--76 (2001).
%
\bibitem{Kra98}
{\sc Kravchenko, O.:} {\it Deformation quantization of symplectic
fibrations.} Compositio Math.\ {\bf 123}, 131--165 (2000).
%
\bibitem{Lic80}
{\sc Lichnerowicz, A.:} {\it Connexions symplectiques et
$\star$-produits invariants.} C.\ R.\ Acad.\ Sc.\ Paris {\bf 291,
A}, 413--417 (1980).
%
\bibitem{Mor86}
{\sc Moreno, C.:} {\it $\ast$-products on some K\"ahler manifolds.}
Lett.\ Math.\ Phys.\ {\bf 11}, 361--372 (1986).
%
\bibitem{MN03a}
{\sc M\"uller, M.\ F., Neumaier, N.:} {\it Some Remarks on
$\mathfrak g$-invariant Fedosov Star Products and Quantum Momentum
Mappings.} Preprint, January 2003, {\bf math.QA/0301101 v2}, to
appear in J.\ Geom.\ Phys..
%
\bibitem{Neu02a}
{\sc Neumaier, N.:} {\it Universality of Fedosov's Construction for
Star Products of Wick Type on Pseudo-K\"ahler Manifolds.} Rep.\
Math.\ Phys.\ {\bf 52}, 43--80 (2003).
%
\bibitem{Sch99}
{\sc Schlichenmaier, M.:} {\it Deformation quantization of compact
K\"ahler manifolds by Berezin-Toeplitz quantization.} in: {\sc
Dito, G., Sternheimer, D.\ (eds.):} Conf\'{e}rence Mosh\'{e} Flato
1999, Vol.\ II.\ Kluwer Academic Publ., Dordrecht, 289--306 (2000).
%
\bibitem{Woo91}
{\sc Woodhouse, N.\ M.\ J.:} {\it Geometric Quantization.} Oxford
Mathematical Monographs.\ Oxford University Press (1991).
%
\bibitem{Xu98}
{\sc Xu, P.:} {\it Fedosov $*$-Products and Quantum Momentum Maps.}
Commun.\ Math.\ Phys.\ {\bf 197}, 167--197 (1998).
\end{thebibliography}
\end{document}